[1994/12/01]% LaTeX date must December 1994 or later
\documentclass[draft]{jotart}
\chardef\bslash=`\\ % p. 424, TeXbook
\usepackage[all]{xy}
\usepackage{amssymb}

\newtheorem{thm}{Theorem}[section]

\newtheorem{lem}[thm]{Lemma}
\newtheorem{prop}[thm]{Proposition}

\newtheorem{example}[thm]{Example}

\theoremstyle{definition}
\newtheorem{defn}{Definition}[section]
\newtheorem{rem}{Remark}[section]

\theoremstyle{remark}

\newcommand{\eval}[2][\right]{\relax
  \ifx#1\right\relax \left.\fi#2#1\rvert}

\begin{document}
\title[C$^*$-algebras of groupoid morphisms]{C$^*$-algebras of groupoid morphisms}
\author[Bruno T. Costa]{Bruno T. Costa}
\email{b.t.costa@ufsc.br}
\author[Renan G. Romano]{Renan G. Romano}
\email{r.g.romano@ufsc.br}
\author[Felipe Vieira]{Felipe Vieira}
\email{f.vieira@ufsc.br}
\address{Departamento de Matem\'{a}tica - UFSC Blumenau\newline
Rua Marechal Rondon, 880. Blumenau, SC. CEP 89065-200. Brazil}

\begin{abstract}
  In this paper we suggest a definition for a C$^*$-algebra attached to an injective morphism of some \'Etale groupoid. We take into account all the peculiarities of such objects and present some interesting relations with already well-known C$^*$-algebras from the literature.
\end{abstract}
\maketitle

\section{INTRODUCTION}

The study of C$^*$-algebras representing distinct mathematical objects is a huge field of research, since all types of the latters can be, somehow, represented as operators over some suitable operator space (\cite{Li1}, \cite{Li2}, \cite{Li3}, \cite{CLY}).

In \cite{Hir}, Hirshberg constructed a $C^*$-algebra associated with endomorphisms of groups with finite cokernel. His inspiration was to see the group elements as unitary operators on a convenient Hilbert space and the endomorphism as one isometry acting on the same object and, among other results, he proved that it can be seen as a crossed product too. In a subsequent paper, Vieira constructed the same $C^*$-algebra for endomorphisms with infinite cokernel (see \cite{Fel}).

Particularly, the use of groupoids and semigroupoids in these processes became convenient, since almost every structure can be viewed as one of them (\cite{Jing}, \cite{Hui}). In this paper, we propose a definition replacing groups by groupoids. The notion of an abstract groupoid was introduced by Brandt in $1927$ (\cite{Bra}) and since then it has been vastly studied. In particular, we are interested in a special kind of groupoid called \textit{\'etale groupoid}, because of its interesting and convenient properties (see \cite{Ren} for more details).

With this picture in mind, we define a $C^*$-algebra associated with an injective morphism $\varphi : G \longrightarrow G$,
of an \'etale groupoid $G$. This case differs to that one of groups because, in a groupoid, we can not multiply every pair of elements. Therefore, the groupoid should be represented in some C$^*$-algebra respecting this condition.

Analyzing a concrete representation of the groupoid and the morphism inside $B(l^2(G))$, it is an easy task to suggest the definition of a universal $C^*$-algebra $\mathbb{U}[\varphi]$, via generators and relations. After, we present and prove some relations between this object and other already known ones.

\section{GROUPOIDS AND SEMIGROUPOIDS}

In order to maximize our intuition, we review the concepts of a groupoid and a semigroupoid as presented in \cite{Mac} and \cite{Exel}. In general, a groupoid $G$ is very similar to a group, as a semigroupoid $S$ is very similar to a semigroup, except that it is not possible to
multiply every two elements of $G$ (or $S$). 
Since some elements are not multipliable in groupoids, its unity is, in fact, a family of elements parametrized by the elements of the base space $M$. A semigroupoid, like a semigroup, does not need a unit. When the semigroupoid has it, we will make it explicit. Formally:
\begin{defn}
A groupoid over a set $M$ is a non-empty set $G$ equipped with a \emph{source projection} $\, \sigma_G^{}: G \longrightarrow M$
and a \emph{target projection} $\, \tau_G^{}: G \longrightarrow M$, a multiplication
\begin{equation*}
 \begin{array}{cccc}
  \cdot: & G^{(2)} & \longrightarrow &  G \\
            &    (h,g)     &   \longmapsto   & h g
 \end{array}
\end{equation*}
defined on the subset
\begin{equation*}
 G^{(2)}~
 =~\bigl\{ (h,g) \in G \times G \,|\, \sigma_G^{}(h)=\tau_G^{}(g) \bigr\}
\end{equation*}
of $G \times G$, such that
\[
\sigma_G(gh)=\sigma_G(h)~~,~~\tau_G(gh)=\tau_G(g)
\]
and satisfying the associativity:
\[
g(hk) = (gh)k \quad \mbox{when} \quad
\sigma_G(g) = \tau_G(h)~~\mbox{e}~~\sigma_G(h)=\tau_G(k) \, ;
\]
a unity
\begin{equation*}
 \begin{array}{cccc}
  1: & M & \longrightarrow &  G  \\[1mm]
     & x &   \longmapsto   & 1_x
 \end{array}
\end{equation*}
such that
\[
\sigma_G(1_x) = x = \tau_G(1_x)
\]
and
\[
g \, 1_{\sigma_G(g)} = g = 1_{\tau_G(g)} \, g \, ;
\]
an inversion
\begin{equation*}
 \begin{array}{cccc}
  \iota_G^{}: & G & \longrightarrow &    G     \\
                & g &  \longmapsto   & g^{-1}
 \end{array}
\end{equation*}
such that
\[
\sigma_G(g^{-1}) = \tau_G(g)~~,~~\tau_G(g^{-1}) = \sigma_G(g)
\]
and
\[
g^{-1} \, g = 1_{\sigma_G(g)}~~,~~g \, g^{-1} = 1_{\tau_G(g)}.
\]
\end{defn}

\begin{defn}
A groupoid $G$ over a base space $M$ is a topological groupoid if both $G$ and $M$ are topological spaces and all structural maps are continuous. 
\end{defn}

\begin{defn}
A topological groupoid is an étale groupoid if the source map is a local homeomorphism. 
\end{defn}

\begin{rem}
In an étale groupoid, since the inverse map is a homeomorphism, the target map is also a local homeomorphism. Even more, for every $x$ in the base space, the source fiber $\sigma^{-1}(x)$ and the target fiber $\tau^{-1}(x)$ are discrete spaces in the subspace topology.
\end{rem}

\begin{rem}
The topology of an étale groupoid $G$ is generated by a base of open sets called open bisections. An open set $U \subseteq G$ is an open bisection if $\sigma$ and $\tau$ restricted to $U$ are homeomorphism onto an open subset of $M$. It means the topology of $G$ behaves locally exactly like the topology of $M$.
\end{rem}

In the next two definitions, we present another key concepts that we need.
\begin{defn}\label{def:groupoidmorphism}
 Let $G$ be a groupoid over $M$ and $H$ a groupoid over $N$ with source projections $\sigma_G$ and $\sigma_H$
 and target projections $\tau_G$ and $\tau_H$, respectively. A \emph{groupoid morphism}, from $G$ to $H$, is a pair of maps
 $\varphi: G \longrightarrow H$ and $\phi: M \longrightarrow N$ satisfying
 $$
 \sigma_H \circ \varphi = \phi \circ \sigma_G \,\, \mbox{ and } \,\, \tau_H \circ \varphi = \phi \circ \tau_G,
 $$
 that is, the diagram
 $$
 \begin{array}{c}
  \xymatrix{
   G ~\ar[r]^{\varphi}
    \ar@<0.5ex>[d]^{\sigma_G} \ar@<-0.5ex>[d]_{\tau_G} & H \ar@<0.5ex>[d]^{\tau_H} \ar@<-0.5ex>[d]_{\sigma_H} \\
   M \ar[r]_{\phi}  &  N
  }
 \end{array}
$$
commutes and $\varphi$ preserves the product: for all $(g,h) \in G^{(2)}$,
$$
\varphi(gh) = \varphi(g)\varphi(h).
$$
Observe that $(\varphi(g),\varphi(h)) \in H^{(2)}$, because
$$
\sigma_H(\varphi(g)) = \phi(\sigma_G(g)) = \phi(\tau_G(h)) = \tau_H(\varphi(h)).
$$
We say that $\varphi$ is a morphism over $\phi$. If $M=N$ and $\phi$ is the identity, we say $\varphi$ is a strict morphism.
\end{defn}

\begin{defn} \label{sgroup}
  A \emph{semigroupoid} is a triple $(\Lambda,\Lambda^{(2)},\cdot)$ such that $\Lambda$ is a set, $\Lambda^{(2)}$ is a subset of $\Lambda\times\Lambda$ and $\cdot:\Lambda^{(2)}\rightarrow\Lambda$ is an operation $(r,s)\mapsto rs$ such that if any of the following three items holds
  \begin{itemize}
    \item[(i)] $(r,s),(s,t)\in\Lambda^{(2)}$;
    \item[(ii)] $(r,s)\in\Lambda^{(2)}$ and $(rs,t)\in\Lambda^{(2)}$;
    \item[(iii)] $(s,t)\in\Lambda^{(2)}$ and $(r,st)\in\Lambda^{(2)}$;
  \end{itemize}
  then $(r,s),(s,t),(rs,t)$ and $(r,st)$ lie in $\Lambda^{(2)}$ and
  \begin{equation*}
    (rs)t=r(st).
  \end{equation*}
\end{defn}

To finish this section, we present three definitions that will be used in the next sections.

\begin{defn}
  Let $r,t\in\Lambda$ where $\Lambda$ is a semigroupoid. We shall say that $r$ \emph{divides} $t$, or that $t$ is a \emph{multiple} of $r$, in symbols $r\,|\,t$, if either
  \begin{itemize}
    \item[(i)] $r=t$ or
    \item[(ii)] there exists $s\in\Lambda$ such that $(r,s)\in\Lambda^{(2)}$ and $rs=t$.
  \end{itemize}
\end{defn}

The item $(i)$ above is necessary since $\Lambda$ may not have a unit.
\begin{defn} We say that $r$ and $s$ in $\Lambda$ \emph{intersect} if they admit a common multiple, that is, there exists $t\in\Lambda$ such that $r|t$ and $s|t$. In this case, we write $r\Cap s$ and, otherwise, we say that $r$ and $s$ are \emph{disjoint} and write $r\perp s$.
\end{defn}

For our last definition, consider $\Lambda^{r}:=\{s\in\Lambda:(r,s)\in\Lambda^{(2)}\}$.
\begin{defn}
  An element $r\in\Lambda$ is \emph{monic} if for every $s,t\in\Lambda^{r}$ we have that $rs=rt$ implies $s=t$ and we say that $r\in\Lambda$ is a \emph{spring} when $\Lambda^{r}=\emptyset$.
\end{defn}

\section{SEMIGROUPOID SEMIDIRECT PRODUCT AND ITS $C^*$-ALGEBRA}\label{secsg}

Starting from an action of the semigroup $\mathbb{N}$ on a groupoid $G$, in this section we show how to construct the semidirect product $G\rtimes_{\varphi}\mathbb{N}$, where $\varphi:G\longrightarrow G$ is an injective groupoid morphism over $\phi: M \longrightarrow M$. We then prove that it is, in fact, a semigroupoid and exhibit its $C^*$-algebra as defined in \cite{Exel}.

\begin{defn} \label{def1}
  Consider $G$ a groupoid and $\varphi:G\rightarrow G$ an injective groupoid morphism over $\phi: M \longrightarrow M$. We define the \emph{semidirect product} $G\rtimes_{\varphi}\mathbb{N}$ as the set $\Lambda=G\times\mathbb{N}$ together with
  $$
      \Lambda^{(2)}:=\{\left((g,m),(h,n)\right)\in\Lambda\times\Lambda:(g,\varphi^{m}(h))\in G^{(2)}\}.
    $$
      The product $\cdot:\Lambda^{(2)}\rightarrow\Lambda$ is defined as follows: given
      $((g,m),(h,n)) \in \Lambda^{(2)}$,
    $$
    (g,m)(h,n):=(g\varphi^{m}(h),m+n).
    $$
\end{defn}

Consider the source projection $\sigma$ and target projection $\tau$ in $G$ and define the functions
$\overline{\sigma}: D_{\overline{\sigma}} \longrightarrow M$ and $\overline{\tau}: D_{\overline{\tau}} \longrightarrow M$ as
  \begin{equation} \label{eq1}
    \begin{split}
      \overline{\sigma}(g,m):&= \phi^{-m}(\sigma(g))\\
      \overline{\tau}(g,m):&=\tau(g).
    \end{split}
  \end{equation}
\begin{prop}\label{pr1}
  The functions $\overline{\sigma}$ and $\overline{\tau}$ have domains
  $$
  D_{\overline{\sigma}}=\{(g,m)\in\Lambda:\sigma(g)\in\phi^{m}(M)\}
  $$
  and $D_{\overline{\tau}}=\Lambda$ respectively. Also, it holds that
  \begin{equation*}
    ((g,m),(h,n))\in\Lambda^{(2)}\Leftrightarrow\,\sigma(g)\in\phi^{m}(M)\textrm{ and }\overline{\sigma}(g,m)=\overline{\tau}(h,n).
  \end{equation*}
\end{prop}

\begin{proof}
We observe that $(g,m) \in D_{\overline{\sigma}}$ if, and only if, there exists $x \in M$ such that $x = \phi^{-m}(\sigma(g))$,
  that is, $\sigma(g) = \phi^m(x) \in \phi^{m}(M)$. It is straightforward that $D_{\overline{\tau}}=\Lambda$.\\
 Besides, note that $((g,m),(h,n))\in\Lambda^{(2)}$ if, and only if, $(g,\varphi^{m}(h))\in G^{(2)}$, that is equivalent to
  $\sigma(g)=\tau(\varphi^{m}(h))=\phi^{m}(\tau (h))$ and so $\sigma(g)\in\phi^{m}(M)$. Now, since $\varphi$ is injective
  -- and then $\phi$ is injective -- we obtain $\phi^{-m}(\sigma(g))=\tau(h)$, that is,
  $\overline{\sigma}(g,m)=\overline{\tau}(h,n)$.
\end{proof}

Now, lets take a look at the set $\Lambda$ and, for this, fix an element $g\in G$ (an arrow that connects $\sigma(g)$ and $\tau(g)$, which are points in the base space). Note first that the arrow $(g,0)$ in the semigroupoid $\Lambda$ can be identified with $g$ since they have the same source and target in $M$. The arrow $(g,1)$ is obtained from $(g,0)$ by changing the source $\sigma(g)$ of $(g,0)$  to $\phi^{-1}(\sigma(g))$, if it is possible, and fixing the target. In general, the arrow $(g,m+n)$ is obtained from $(g,m)$ by applying $\phi^{-n}$ on the source of $(g,m)$.

On the other hand, the element $(\varphi^{n}(g),0)$ is obtained from $(g,0)$ by applying $\phi^{n}$ on both source and target of
$(g,0)$, since $\varphi$ is a morphism. Finally, the element $(\varphi^{n}(g),m)$ for $m\leq n$ is obtained from $(g,0)$ by applying $\phi^{n}$ on the target and $\phi^{n-m}$ on the source of $(g,0)$.

\begin{prop} \label{pr2}
  If $((g,m),(h,n))\in\Lambda^{(2)}$ and $(h,n)\in D_{\overline{\sigma}}$ then $
  (g,m)(h,n)\in D_{\overline{\sigma}}$, $\overline{\sigma}((g,m)(h,n))=\overline{\sigma}(h,n)$ and $\overline{\tau}((g,m)(h,n))=\overline{\tau}(g,m)$.
\end{prop}

\begin{proof}
  Note that if $(h,n)\in D_{\overline{\sigma}}$ then $\sigma(h)\in\phi^{n}(M)$. This implies that
  $$
  \phi^{m}(\sigma(h))\in\phi^{m+n}(M).
  $$
  Using the properties of $\sigma$ and the fact that $\varphi$ is a groupoid morphism, we obtain
  \begin{equation*}
    \phi^{m+n}(M)\ni\phi^{m}(\sigma (h))=\sigma(\varphi^{m}(h))=\sigma(g\varphi^{m}(h))
  \end{equation*}
  since $(g,\varphi^{m}(h))\in G^{(2)}$. This implies that $(g,m)(h,n)=(g\varphi^{m}(h),m+n)\in D_{\overline{\sigma}}$. Now,
  \begin{equation*}
  \begin{split}
    \overline{\sigma}((g,m)(h,n))&=\phi^{-(m+n)}(\sigma(g\varphi^{m}(h)))\\
    &=\phi^{-(m+n)}(\sigma(\varphi^{m}(h)))\\
    &=\phi^{-n}(\sigma(h))\\
    &=\overline{\sigma}(h,n).
  \end{split}
  \end{equation*}
  For $\overline{\tau}$ we have
  \begin{equation*}
    \overline{\tau}((g,m)(h,n))=\overline{\tau}(g\varphi^{m}(h),m+n)=\tau(g\varphi^{m}(h))=\tau(g)=\overline{\tau}(g,m).
  \end{equation*}
\end{proof}

The next lemma will be important to show that $\Lambda$, as constructed here, is a semigroupoid.
\begin{lem} \label{lema1}
 If $((g,m),(h,n))\in\Lambda^{(2)}$ and $(g,m)(h,n)\in D_{\overline{\sigma}}$, then $(h,n)\in D_{\overline{\sigma}}$.
\end{lem}

\begin{proof}
 We have, using the injectivity of $\phi$, that
 \begin{equation*}\begin{split}
 (g,m)(h,n)\in D_{\overline{\sigma}} &\Rightarrow (g\varphi^m(h), m+n) \in D_{\overline{\sigma}}\\
 &\Rightarrow \sigma(g\varphi^m(h)) \in \phi^{m+n}(M) \\
                                     &\Rightarrow \sigma(\varphi^m(h)) \in \phi^{m+n}(M)\\
 &\Rightarrow \phi^m(\sigma(h)) \in \phi^{m+n}(M) \\
                                     &\Rightarrow \sigma(h) \in \phi^n(M)\\
 &\Rightarrow (h,n) \in D_{\overline{\sigma}}.
 \end{split}\end{equation*}
\end{proof}

\begin{thm}
  The triple $(\Lambda,\Lambda^{(2)},\cdot)$ is a semigroupoid.
\end{thm}

\begin{proof}
  By Proposition \ref{pr2} and Lemma \ref{lema1}, given $((f,l),(g,m))$,\linebreak
  $((g,m),(h,n))$, $((f,l)(g,m),(h,n))$
  and $((f,l),(g,m)(h,n))$ in $\Lambda \times \Lambda$, if any of the items of the Definition \ref{sgroup}
  hold, then $((f,l),(g,m))$, $((g,m),(h,n))$,\linebreak $((f,l)(g,m),(h,n))$ and $((f,l),(g,m)(h,n))$ lie in $\Lambda^{(2)}$.
  Moreover, since $G$ is a groupoid,
    \begin{equation*}
      \begin{split}
        [(f,l)(g,m)](h,n)&=(f\varphi^{l}(g),l+m)(h,n)\\
        &=(f\varphi^{l}(g)\varphi^{l+m}(h),l+m+n)\\
        &=(f\varphi^{l}(g\varphi^{m}(h)),l+m+n)\\
        &=(f,l)(g\varphi^{m}(h),m+n)\\
        &=(f,l)[(g,m)(h,n)].
      \end{split}
    \end{equation*}
\end{proof}

In the following, we present some interesting properties of $\Lambda$, since non monic elements and springs can not be properly represented in a C$^*$-algebra.
\begin{prop}
  Every element of $\Lambda$ is monic.
\end{prop}

\begin{proof}
  Suppose that $(g,m),(h,n)\in\Lambda^{(f,l)}$ with $(f,l)(g,m)=(f,l)(h,n)$. Then $(f\varphi^{l}(g),l+m)=(f\varphi^{l}(h),l+n)$, that is, $m=n$ and $f\varphi^{l}(g)=f\varphi^{l}(h)$. Using $f^{-1}$ and the injectivity of $\varphi$, we conclude that $g=h$ and $(f,l)$ is monic.
\end{proof}

\begin{prop}
  An element $(g,m)\in\Lambda$ is a spring if, and only if, $$(g,m)\notin D_{\overline{\sigma}}.$$
\end{prop}

\begin{proof}
  Observe that $(g,m)$ is a spring if, and only if, $\sigma(g)\neq\phi^{m}(\tau(h))$ for all $h\in G$, which means that $(g,m)\notin D_{\overline{\sigma}}$.
\end{proof}

Therefore we will consider $D_{\overline{\sigma}}=\Lambda$ because, in this setup, the semigroupoid $\Lambda$ will have no springs. This is equivalent to say that the function $\phi$ in the Definition \ref{def:groupoidmorphism} is bijective.
\begin{example}
 Given a set $M$, consider the cartesian product $M \times M$ and define the source projection $\sigma(y,x) := x$ and target projection $\tau(y,x) := y$ The multiplication will be $(z,y)(y,x) := (z,x)$, the diagonal as unit, $1_x := (x,x)$, and switch as inversion, $(y,x)^{-1} := (x,y)$. Then $M \times M$ is a groupoid over~$M$, called the \emph{pair groupoid} of~$M$. Given a bijective map $\phi:M \longrightarrow M$ we can construct a groupoid morphism $\varphi: M \times M \longrightarrow M \times M$ given by:
 $$
 \varphi(y,x) := (\phi(y),\phi(x))
 $$
 where the following diagram commutes
 $$
 \begin{array}{c}
  \xymatrix{
   M \times M ~\ar[r]^{\varphi}
    \ar@<0.5ex>[d]^\sigma \ar@<-0.5ex>[d]_\tau & M \times M \ar@<0.5ex>[d]^\tau \ar@<-0.5ex>[d]_\sigma \\
   M \ar[r]_{\phi}  &  M
  }
 \end{array}
$$
Observe that, since $\phi$ is injective, $\varphi$ is injective as well.\\[0.2cm]
This leads us to the construction of a semigroupoid $(M \times M) \rtimes_\varphi \mathbb{N}$, the semidirect product
of $M \times M$ with $\mathbb{N}$. Namely, using the maps defined in (\ref{eq1}) and multiplication from Definition \ref{def1} we obtain, for
$(y,x,n) \in (M \times M) \rtimes_\varphi \mathbb{N}$, a source projection
$$
\bar{\sigma}(y,x,n) = \phi^{-n}(\sigma(y,x)) = \phi^{-n}(x),
$$
a target projection
$$
\bar{\tau}(y,x,n) = \tau(y,x) = y
$$
and multiplication
$$
\begin{aligned}
(y_2,x_2,n_2)(y_1,x_1,n_1) &= ((y_2,x_2)\varphi^{n_2}(y_1,x_1), n_1 + n_2) \\
                                  &=((y_2,x_2)(\phi^{n_2}(y_1),\phi^{n_2}(x_1)), n_1 + n_2) \\
                                  &= (y_2,\phi^{n_2}(x_1),n_1 + n_2),
\end{aligned}
$$
where $y_1 = \phi^{-n_2}(x_2)$.
The set $(M \times M) \rtimes_\varphi \mathbb{N}$ endowed with this structure is a semigroupoid.
\end{example}

Since we want to present an isomorphism involving the C$^*$-algebra of a semigroupoid, we repeat its construction (from \cite{Exel}) in the following two definitions.
\begin{defn}\label{defi0}
  Let $\Lambda$ be a semigroupoid and $B$ a unital $C^{*}$-algebra. A map $S:\Lambda\rightarrow B$ given by $r\mapsto S_{r}$ will be called a \emph{representation} of $\Lambda$ in $B$ if, for every $r,t\in\Lambda$, denoting by $P_{r}=S_{r}S_{r}^{*}$ and $Q_{r}=S_{r}^{*}S_{r}$ the initial and final projections, respectively, we have:
  \begin{itemize}
    \item[(i)] $S_{r}$ is a partial isometry;
    \item[(ii)] $S_{r}S_{t}=\begin{cases}
      S_{rt}\textrm{ if }(r,t)\in\Lambda^{(2)}\\
      0,\textrm{ otherwise;}
    \end{cases}$
    \item[(iii)] $P_{r}P_{t} = 0$ if $r\perp t$;
    \item[(iv)] $Q_{r}P_{t}=P_{t}$ if $(r,t)\in\Lambda^{(2)}$.
  \end{itemize}
\end{defn}

In order to define the $C^*$-algebra of a semigroupoid, normally one considers tight representations. However, we have observed that, in \cite{Exel}, the main reason to choose tight representations (and not merely a representation) is to ensure an isomorphism between the $C^*$-algebra of the Markov semigroupoid and the unital Cuntz-Krieger algebra introduced in \cite{ExLa}. Therefore, we will not use tight representations to define the $C^*$-algebra of a semigroupoid, since it is not needed to derive our results.
\begin{defn}
Given a semigroupoid $\Lambda$, its $C^*$-algebra $C^*(\Lambda)$ is the universal unital C$^*$-algebra generated by a family of partial
isometries $\{S_r\}_{r \in \Lambda}$ such that $S$ is a representation of $\Lambda$.
\end{defn}

In the next section we would like to prove there is an isomorphism between the $C^*$-algebra of a semigroupoid and the C$^*$-algebra of a groupoid endomorphism. However, on the one hand the C$^*$-algebra of the groupoid endomorphism is naturally equipped with an isometry and on the other hand the $C^*$-algebra of the semigroupoid does not explicitly contain an isometry that can be used to construct such an isomorphism. At least in the case where the base space $M$ is not finite. To overcome this difficulty, we define the following universal C$^*$-algebra: 

\begin{defn}
Consider $G$ a groupoid and $\varphi:G\rightarrow G$ an injective groupoid morphism over $\phi: M \longrightarrow M$. Given the semigroupoid $\Lambda = G\rtimes_{\varphi}\mathbb{N}$ and an isometry $s$, the $C^*$-algebra $C^*(\Lambda,s)$ is the universal unital C$^*$-algebra generated by a family of partial isometries $\{S_r\}_{r \in \Lambda}$ such that $S$ is a representation of $\Lambda$, and an isometry $s$ satisfying
\begin{equation}\label{isos}
sS_{(g,m)} = S_{(\varphi(g),m)}s.
\end{equation}
\end{defn}

\section{THE C$^*$-ALGEBRA OF A GROUPOID MORPHISM}

Based on the similarities of groupoids and groups, we adapt the definition of the C$^*$-algebra of a group endomorphism to our setup. Let $G$ be an étale groupoid over $M$ and $G_x = \sigma^{-1}(x)$ be the discrete source fiber. In addition, consider the Hilbert space $\ell^2(G_x)$ with counting measure and its canonical orthonormal basis $\{\xi_h:h \in G_x\}$. Observe that this basis can be constructed since $G$ is étale, that is, $G_x = \sigma^{-1}(x)$ is discrete. Let $\varphi:G\rightarrow G$ be a continuous injective groupoid morphism. We define the Hilbert space $\bigoplus_{x \in M} \ell^2(G_x)$. For each $g\in G$, define the operator $U_g$: 
$$
U_g(\xi_h) =
\begin{cases}
      \xi_{gh}, & \mbox{if } (g,h) \in G^{(2)} \\
\,\,       0,   & \mbox{otherwise.}
\end{cases}
$$
It is easy to prove that they are partial isometries. Besides, the endomorphism $\varphi$ induces an isometry $S$
acting on $\bigoplus_{x \in M} \ell^2(G_x)$ by
$$
S(\xi_h) = \xi_{\varphi(h)}.
$$
The reduced C$^{*}$-algebra associated with $\varphi$ is the C$^{*}$-subalgebra of $$\mathcal{B}\left(\bigoplus_{x \in M} \ell^2(G_x)\right)$$ generated by $U_g$ and $S$. its elements satisfy interesting properties which are the base of our following definition.

\begin{defn}\label{defi1} Consider an injective endomorphism $\varphi$ of an \'etale groupoid $G$. We define $\mathbb{U}[\varphi]$ as the universal $C^*$-algebra generated by
$$
\mbox{partial isometries }\{u_g\,:\,g\in G\} \mbox{ and}
$$
$$
\mbox{an isometry }s
$$
satisfying
\begin{enumerate}
\item[(i)]
$u_gu_h =
\begin{cases}
u_{gh}, & \mbox{if } (g,h) \in G^{(2)} \\
\,\, 0, & \mbox{otherwise};
\end{cases}
$\\
\item[(ii)] $u_g^*=u_{g^{-1}}$;\\
\item[(iii)] $su_g=u_{\varphi(g)}s$;\\
%\item[(iv)] $\displaystyle\sum_{x\in M} u_{1_{x}}s=s$.
\end{enumerate}
Furthermore, the elements $u_gs^m{s^*}^mu_{g^{-1}}$ and $u_hs^n{s^*}^nu_{h^{-1}}$ are projections that commute and
if $m \leqslant n$ then
\begin{enumerate}
\item[(iv)] $u_gs^m{s^*}^mu_{g^{-1}}u_hs^n{s^*}^nu_{h^{-1}}=0$, if $g^{-1}h \notin \varphi^m(G)$.
\end{enumerate}
\end{defn}

In (iii), observe that $u_{g^{-1}}s^*=s^*u_{\varphi(g^{-1})}$.

Along this section and the next, we will consider $\varphi$ is bijective. Our goal now is to show that the $C^*$-algebra $\mathbb{U}[\varphi]$ constructed in this section is isomorphic to the $C^*$-algebra $C^*(\Lambda,s)$, presented in Section \ref{secsg}. For this, we need an easy lemma.
\begin{lem} \label{INT}
It holds that $g^{-1}h \in \varphi^m(G)$ if, and only if, $(g,m) \Cap (h,n)$.
\end{lem}
\begin{proof}
If $g^{-1}h \in \varphi^m(G)$ then there exists $f \in G$ such that $g^{-1}h = \varphi^m(f)$, i.e., $h = g\varphi^m(f)$. Also, note that there exist $l,k \in \mathbb{N}$ such that $m + l = n + k$. So,
\begin{equation*}\begin{split}
(g,m)(f,l) &= (g\varphi^m(f),m+l) \\
&= (h,m+l) \\
&= (h1_{\sigma(h)},m+l) \\
&= (h1_{\phi^n(x)},m+l) \\
&= (h\varphi^n(1_x),m+l) \\
&= (h\varphi^n(1_x),n+k) \\
&= (h,n)(1_x,k).
\end{split}\end{equation*}
It means that there exist $(f,l),(1_x,k) \in \Lambda$ such that $(g,m)(f,l) = (h,n)(1_x,k)$, i.e., $(g,m) \Cap (h,n)$. Now, if
$(g,m) \Cap (h,n)$ it is obvious that $g^{-1}h \in \varphi^m(G)$.
\end{proof}

\begin{thm} \label{TISO}
$\mathbb{U}[\varphi]$ is isomorphic to $C^*(\Lambda,s)$.
\end{thm}
\begin{proof}
To start, let us use the notation of Definition \ref{defi0} to define
$$
\Phi: \mathbb{U}[\varphi] \longrightarrow C^*(\Lambda,s),
$$
given by $\Phi(u_g) := S_{(g,0)}$ and $\Phi(s) := s$. We need to show that these representatives satisfy the relations generating $\mathbb{U}[\varphi]$.\\[0.2cm]
Using the condition (ii) of Definition \ref{defi0}, relation (i) generating $\mathbb{U}[\varphi]$ follows by easy calculations.\\[0.2cm]
Condition $(ii)$ follows from condition $(i)$ of Definition \ref{defi0}: on one hand, since $S_{(g,0)}$ is a partial isometry,
$S_{(g,0)}S^*_{(g,0)}S_{(g,0)}=S_{(g,0)}$. On the other hand,
\begin{equation*}\begin{split}
S_{(g,0)}S_{(g^{-1},0)}S_{(g,0)} &= S_{(g,0)(g^{-1},0)}S_{(g,0)}\\
&= S_{(gg^{-1},0)}S_{(g,0)}\\
&= S_{(1_{\tau(g)},0)}S_{(g,0)}\\
&= S_{(g,0)}.
\end{split}\end{equation*}
Therefore, $S^*_{(g,0)}=S_{(g^{-1},0)}$.\\
Finally, (iii) is straightforward from (\ref{isos})  and (iv) follows from Lemma \ref{INT} and condition
(iii) of the Definition \ref{defi0}. \\[0.2cm]
For the inverse map, consider
$$
\Psi: C^*(\Lambda,s) \longrightarrow \mathbb{U}[\varphi],
$$
given by $\Psi(S_{(g,m)}) := u_gs^m$ and $\Psi(s) = s$. Let us show that the elements $u_gs^m$ and $s$ satisfy the five relations that generate $C^*(\Lambda,s)$.\\
To prove condition (i) note that
\begin{equation*}\begin{split}
(u_gs^m)(u_gs^m)^*(u_gs^m) &= u_gs^m{s^*}^mu_{g^{-1}}u_gs^m\\
& = u_gs^m{s^*}^mu_{1_{\sigma(g)}}s^m\\
& = u_gs^mu_{\varphi^{-m}(1_{\sigma(g)})}{s^*}^ms^m \\
&= u_gs^mu_{\varphi^{-m}(1_{\sigma(g)})}\\
& = u_gu_{1_{\sigma(g)}}s^m\\
& = u_gs^m.
\end{split}\end{equation*}
Likewise, we can prove that
$$
(u_gs^m)^*(u_gs^m)(u_gs^m)^*=(u_gs^m)^*.
$$
The condition (ii) follows from (i) of Definition \ref{defi1}.\\
For (iii), if $\tau(g)\neq\tau(h)$ then
$$
(u_gs^m{s^*}^mu_{g^{-1}})(u_hs^n{s^*}^nu_{h^{-1}})=0,
$$
as we want to proof. Now, if $\tau(g)=\tau(h)$ and $(g,m)\perp (h,n)$, Lemma \ref{INT} implies that $g^{-1}h \notin \varphi^m(G)$. In this case,
\begin{equation*}
u_gs^m{s^*}^mu_{g^{-1}}u_hs^n{s^*}^nu_{h^{-1}} = 0.
\end{equation*}
To show (iv) consider $((g,m),(h,n)) \in \Lambda^{(2)}$. Then,
\begin{equation*}\begin{split}
({s^*}^mu_{g^{-1}}u_gs^{m})(u_hs^{n}{s^*}^nu_{h^{-1}})&={s^*}^mu_{g^{-1}}u_{g}u_{\varphi^m(h)}s^{m}s^{n}{s^*}^nu_{h^{-1}}\\
&= {s^*}^mu_{\varphi^m(h)}s^{m+n}{s^*}^nu_{h^{-1}}\\
&= {s^*}^ms^mu_hs^n{s^*}^nu_{h^{-1}}\\
&= u_hs^n{s^*}^nu_{h^{-1}}.
\end{split}\end{equation*}
The condition (\ref{isos}) follows from (iii) of Definition \ref{defi1}:
$$
su_gs^m=u_{\varphi(g)}s^{m+1} = u_{\varphi(g)}s^ms.
$$
Finally, a simple calculation shows that $\Psi$ is the inverse map of $\Phi$:
\begin{equation*}\begin{split}
\Psi(\Phi(u_g)) &= \Psi(S_{(g,0)}) = u_g,\\
\Psi(\Phi(s)) &= \Psi(s) = s.
\end{split}\end{equation*}
Besides, an induction argument implies that
$$
\Phi(\Psi(S_{(g,m)})) = S_{(g,m)} \mbox{\,\, and \,\,} \Phi(\Psi(s)) = s
$$
\end{proof}

Observe that we do not impose any additional hypotesis on the base space $M$. We only assume that $G$ is an étale groupoid over $M$. However, if we consider $M$ as a finite topological space, then we are able to construct the isometry $s$ in the original $C^*(\Lambda)$ as
$$
s = \displaystyle\sum_{x\in M} S_{(1_x,1)}.
$$
which converges since $M$ is finite. In fact, $s$ is an isometry: by conditions (i) and (iii) of Definition \ref{defi0}, if $x \neq y$, 
$$
S^*_{(1_x,1)}S_{(1_y,1)} = S^*_{(1_x,1)}P_{(1_x,1)}P_{(1_y,1)}S_{(1_y,1)} = 0.
$$
Now, observe that $((1_x,1),(1_{\phi^{-1}(x)},0)) \in \Lambda^{(2)}$ Then, by condition (ii) of Definition \ref{defi0},
$$
S_{(1_x,1)}S_{(1_{\phi^{-1}(x)},0)} = S_{(1_x,1)(1_{\phi^{-1}(x)},0)} = S_{(1_x,1)}.
$$
Furthermore, since $S_{(1_{\phi^{-1}(x)},0)}$ is a projection,
$$
P_{(1_{\phi^{-1}(x)},0)} = S_{(1_{\phi^{-1}(x)},0)}.
$$
Then, by condition (iv)
$$
Q_{(1_x,1)}S_{(1_{\phi^{-1}(x)},0)} = S_{(1_{\phi^{-1}(x)},0)}.
$$
Considering that $S_{(1_{\phi^{-1}(x)},0)}$ is a projection, we obtain
$$
S_{(1_x,1)} = S_{(1_x,1)}P_{(1_{\phi^{-1}(x)},0)},
$$
that is,
$$
S^*_{(1_x,1)} = P_{(1_{\phi^{-1}(x)},0)}S^*_{(1_x,1)}.
$$
Thus,
$$
\begin{aligned}
Q_{(1_x,1)} = S^*_{(1_x,1)}S_{(1_x,1)} &= 
P_{(1_{\phi^{-1}(x)},0)}S^*_{(1_x,1)}S_{(1_x,1)}P_{(1_{\phi^{-1}(x)},0)} \\
&= P_{(1_{\phi^{-1}(x)},0)}Q_{(1_x,1)}P_{(1_{\phi^{-1}(x)},0)} \\
&= S_{(1_{\phi^{-1}(x)},0)}Q_{(1_x,1)}S_{(1_{\phi^{-1}(x)},0)}.
\end{aligned}
$$
Finally,
$$
Q_{(1_x,1)} = S_{(1_{\phi^{-1}(x)},0)}Q_{(1_x,1)}S_{(1_{\phi^{-1}(x)},0)}
= S_{(1_{\phi^{-1}(x)},0)}S_{(1_{\phi^{-1}(x)},0)} = S_{(1_{\phi^{-1}(x)},0)}.
$$
We are able to prove that $s^*s$ is the identity of $C^*(\Lambda)$:
$$
s^*s = \displaystyle\sum_{x\in M}\displaystyle\sum_{y\in M} S^*_{(1_x,1)}S_{(1_y,1)}
= \displaystyle\sum_{x\in M} S^*_{(1_x,1)}S_{(1_x,1)}
= \displaystyle\sum_{x\in M} Q_{(1_x,1)} 
= \displaystyle\sum_{x\in M} S_{(1_{\phi^{-1}(x)},0)}
$$
which is clearly the identity of $C^*(\Lambda)$. It proves that $s$ is an isometry.

It remains to verify that
$$
sS_{(g,m)} = S_{(\varphi(g),m)}s.
$$
Observe that, on the one hand
$$
sS_{(g,m)} = \displaystyle\sum_{x\in M} S_{(1_x,1)}S_{(g,m)}
= S_{(\varphi(g),m+1)}.
$$
On the other hand,
$$
S_{(\varphi(g),m)}s = \displaystyle\sum_{x\in M} S_{(\varphi(g),m)}S_{(1_x,1)}
= S_{(\varphi(g),m+1)}.
$$
This ensures that $C^*(\Lambda)$ is isomorphic to $\mathbb{U}[\varphi]$ when $M$ is finite.

\section{CROSSED PRODUCT DESCRIPTION}

In this section we want to show that $\mathbb{U}[\varphi]$ can be viewed as a quotient of a crossed product. Due to the particularities of our object, the crossed product that we will construct is a semigroup one, in which the initial C$^*$-algebra is a groupoid crossed product.

This last construction proceeds as follows: given a $C^{*}$-algebra $A$ and a groupoid $G$ over a locally compact Hausdorff space $M$, we can construct a\linebreak groupoid crossed product $A\rtimes_{\alpha}G$ following the steps on \cite{GG}. Basically, because  of the fibered nature of a groupoid, $A$ must have a ``fibration'' so that we can define an action $\alpha$ of $G$ on $A$.

Now, by a $C_{0}(M)$-algebra we mean a $C^{*}$-algebra $A$ equipped with a nondegenerate $*$-homomorphism
$$
C_{0}(M) \longrightarrow ZM(A),
$$
where $ZM(A)$ denotes the center of the multiplier algebra of $A$. Associated with every $C_{0}(M)$-algebra $A$, there is an upper-semicontinuous $C^{*}$-bundle $\mathfrak A$ over $M$, whose fibers are given by
$$
A_x=A/I_x,\qquad
I_x=\overline{C_0(M\setminus{x})A}.
$$
By the representation theorem for $C_{0}(M)$-algebras (see Theorem C.26 in \cite{DPW}), the $C^{*}$-algebra $A$ is canonically isomorphic to the algebra of continuous sections of $\mathfrak A$ vanishing at infinity:
$$
A\cong\Gamma_0(M,\mathfrak A).
$$

The groupoid dynamical system $(A,G,\alpha)$ is defined by choosing an action $\alpha$ of $G$ on $\mathfrak{A}$ and the groupoid crossed product $A\rtimes_{\alpha}G$ is defined as the completion of $\Gamma_{c}(G,\tau^{*}\mathfrak{A})$ with respect to an especific norm.

In our present case, we define the $C^{*}$-subalgebra of $\mathbb{U}[\varphi]$
\begin{equation*}
\mathfrak{D}[\varphi] := C^*(\{u_hs^n{s^*}^nu_{h^{-1}} : h \in G, n \in \mathbb{N}\}).
\end{equation*}
For each $x \in M$, let
\begin{equation*}
\mathfrak{D}_x[\varphi] := C^*(\{u_hs^n{s^*}^nu_{h^{-1}} : h \in G, n \in \mathbb{N} \mbox{ and }\tau(h) = x\})
\end{equation*}
which gives us the $C^{*}$-bundle
\begin{equation*}
\mathfrak{D}[\varphi] := \dot{\bigcup_{x\in M}}\mathfrak{D}_x[\varphi]
\end{equation*}
whose projection $\pi: \mathfrak{D}[\varphi] \longrightarrow M$ is defined by $\pi(u_hs^n{s^*}^nu_{h^{-1}}) := \tau(h)$.

Following \cite{GG}, the notation is $\mathfrak{D}[\varphi]=\mathfrak{A}$ and $A$ is $\Gamma_{0}(M,\mathfrak{D}[\varphi])$, which we will denote by $D[\varphi]$ to simplify. Here, the pull-back of $\mathfrak{D}[\varphi]$ by $\tau$ is given by
 $$
 \tau^*\mathfrak{D}[\varphi] := \{(P,g) \in \mathfrak{D}[\varphi] \times G ;\, \tau(g) = \pi(P)\}
 $$
represented by the commutative diagram
 $$
  \begin{array}{c}
   \xymatrix{
    ~\tau^*\mathfrak{D}[\varphi] ~\ar[r]
     \ar[d] & ~\mathfrak{D}[\varphi] ~\ar[d]^\pi \\
    G \ar[r]_\tau  &  M
   }
  \end{array}
 $$
 It remains to be proven that $D[\varphi]$ is a $C_{0}(M)$-algebra. This ensures $\pi: \mathfrak{D}[\varphi] \longrightarrow M$ is a upper-semicontinuous $C^{*}$-bundle, by Theorem C.26 in \cite{DPW}. For that, consider 
 $$
 \Phi: C_0(M) \longrightarrow ZM(D[\varphi])
 $$
given, for each $s \in D[\varphi]$ and $x \in M$, by:
$$
(\Phi(f)s)(x) = f(x)s(x).
$$
Since it is defined as a scalar multiplication, $\Phi$ is a homomorphism. It remains to prove that \(\Phi\) is nondegenerate. Let \(s\in D[\varphi]\) and \(\varepsilon>0\). Since $s$ vanishes at infinity, there exists a compact set \(K\subseteq M\) such that
\[
\|s(x)\|<\varepsilon
\qquad\text{for all }x\in M\setminus K.
\]
Choose \(f\in C_c(M)\) such that \(0\leq f\leq 1\) and \(f(x)=1\) for every \(x\in K\). Then
\[
\|\Phi(f)s-s\| = \sup_{x\in M}|f(x)-1|\,\|s(x)\| \leq \varepsilon.
\]
Consequently,
\[
\overline{\Phi(C_0(M))D[\varphi]}=D[\varphi],
\]
so \(\Phi\) is nondegenerate. Therefore, \(D[\varphi]\) is a \(C_0(M)\)-algebra.

Now, we define the action $\alpha$ from $G$ on the fibred $C^{*}$-algebra $\mathfrak{D}[\varphi]$ as the family of functions $\{\alpha_g\}_{g\in G}$ such that
\begin{equation*}
\alpha_g(u_hs^n{s^*}^nu_{h^{-1}}) := u_gu_hs^n{s^*}^nu_{h^{-1}}u_{g^{-1}} = u_{gh}s^n{s^*}^nu_{(gh)^{-1}}.
\end{equation*}

\begin{prop} Let $\alpha$ be defined as previsouly. Then:
\begin{itemize}
\item[(i)] For all $(f,g) \in G^{(2)}$, we have $\alpha_{fg} = \alpha_f \circ \alpha_g$;
\item[(ii)] Given $f \in G$, $\alpha_f: \mathfrak{D}_{\sigma(f)}[\varphi] \longrightarrow \mathfrak{D}_{\tau(f)}[\varphi]$ is an isomorphism.
\end{itemize}
That is, $\alpha$ is an action from $G$ on the $C^{*}$-algebra $D[\varphi]$.
\end{prop}
\begin{proof} $(i)$: Let $u_hs^n{s^*}^nu_{h^{-1}}\in\mathfrak{D}[\varphi]$. Then,
\begin{equation*}\begin{split}
\alpha_{fg}(u_hs^n{s^*}^nu_{h^{-1}})&=u_{fgh}s^n{s^*}^nu_{(fgh)^{-1}}\\
&=u_f(u_g(u_{h}s^n{s^*}^nu_{h^{-1}})u_{g^{-1}})u_{f^{-1}}\\
&=\alpha_f\circ\alpha_g(u_hs^n{s^*}^nu_{h^{-1}}).
\end{split}\end{equation*}
$(ii)$ Item $(i)$ already ensures that $\alpha_f$ is bijective with inverse $\alpha_{f^{-1}}$. To prove that it is a morphism, consider $u_gs^m{s^*}^mu_{g^{-1}}$ and $u_hs^n{s^*}^nu_{h^{-1}}$ in $\mathfrak{D}_{\sigma(f)}[\varphi]$. Then,
\begin{equation*}\begin{split}
\alpha_{f}((u_gs^m{s^*}^mu_{g^{-1}})&(u_hs^n{s^*}^nu_{h^{-1}}))\\
&=u_f(u_gs^m{s^*}^mu_{g^{-1}}u_hs^n{s^*}^nu_{h^{-1}})u_{f^{-1}}\\
&=u_fu_gs^m{s^*}^mu_{g^{-1}}u_{f^{-1}}u_fu_hs^n{s^*}^nu_{h^{-1}}u_{f^{-1}}\\
&=(u_fu_gs^m{s^*}^mu_{g^{-1}}u_{f^{-1}})(u_fu_hs^n{s^*}^nu_{h^{-1}}u_{f^{-1}})\\
&=\alpha_{f}(u_gs^m{s^*}^mu_{g^{-1}})\alpha_{f}(u_hs^n{s^*}^nu_{h^{-1}}).
\end{split}\end{equation*}
Also
\begin{equation*}\begin{split}
\alpha_{f}((u_gs^m{s^*}^mu_{g^{-1}})^{*})&=\alpha_{f}(u_gs^m{s^*}^mu_{g^{-1}})\\
&=u_f(u_gs^m{s^*}^mu_{g^{-1}})u_{f^{-1}}\\
&=(u_f(u_gs^m{s^*}^mu_{g^{-1}})u_{f^{-1}})^{*}\\
&=(\alpha_{f}(u_gs^m{s^*}^mu_{g^{-1}}))^{*}.
\end{split}\end{equation*}
\end{proof}

With this, we are able to present the groupoid dynamical system $(D[\varphi],G,\alpha)$ and, finally, to construct the groupoid
crossed product $D[\varphi] \rtimes_{\alpha} G$ as the universal enveloping algebra of $\Gamma_{c}(G,\tau^*\mathfrak{D}[\varphi])$.

The next step is to construct the crossed product of $D[\varphi] \rtimes_\alpha G$ with $\mathbb{N}$. Given $n\in\mathbb{N}$, consider $\beta_n: D[\varphi] \rtimes_\alpha G \longrightarrow D[\varphi] \rtimes_\alpha G$ defined by
$$
\beta_n(u_gs^m{s^*}^mu_{g^{-1}},h) := (s^nu_gs^m{s^*}^mu_{g^{-1}}{s^*}^n,\varphi^n(h)).
$$
Then, we construct the crossed product $(D[\varphi] \rtimes_\alpha G) \rtimes_\beta \mathbb{N}$, using the definitions of \cite{Li1}.

\begin{thm}
The $C^*$-algebra $\mathbb{U}[\varphi]$ is a quotient of $(D[\varphi] \rtimes_{\alpha} G) \rtimes_\beta \mathbb{N}$.
\end{thm}
\begin{proof}
The definition of semigroup crossed products ensures the existence of a unital $*$-homomorphism
$$
i_D : D[\varphi]\rtimes_{\alpha} G \longrightarrow (D[\varphi] \rtimes_{\alpha} G) \rtimes_\beta \mathbb{N}
$$
and a homomorphism of semigroups
$$
i_{N} : \mathbb{N} \longrightarrow Isom((D[\varphi] \rtimes_{\alpha} G) \rtimes_\beta \mathbb{N})
$$
satisfying
$$
i_D(\beta_{n}(x)) = i_{N}(n)i_D(x)(i_{N}(n))^*.
$$
But, if we consider the maps
$$
\begin{array}{cccc}
  \varpi:
  & D[\varphi]\rtimes_{\alpha} G  & \longrightarrow & \mathbb{U}[\varphi]    \\[1mm]
  &    (u_gs^m{s^*}^mu_{g^{-1}},h)    &   \longmapsto   &(u_gs^m{s^*}^mu_{g^{-1}})u_h
 \end{array}
$$
and
$$
\begin{array}{cccc}
  \rho:
  & \mathbb{N}   & \longrightarrow & Isom(\mathbb{U}[\varphi]),  \\[1mm]
  &     n      &   \longmapsto    &  s^n
 \end{array}
$$
we observe that
\begin{equation*}\begin{split}
\rho(n)\varpi(u_gs^m{s^*}^mu_{g^{-1}},h)(\rho(n))^*&=s^n((u_gs^m{s^*}^mu_{g^{-1}})u_h)(s^n)^*\\
&=s^nu_gs^m{s^*}^mu_{g^{-1}}{s^*}^n u_{\varphi^n(h)}\\
&=\varpi(s^nu_gs^m{s^*}^mu_{g^{-1}}{s^*}^n,\varphi^n(h))\\
&=\varpi(\beta_{n}(u_gs^m{s^*}^mu_{g^{-1}},h)),
\end{split}\end{equation*}
that is, $(\mathbb{U}[\varphi],\varpi,\rho)$ is a covariant representation of $((D[\varphi] \rtimes_{\alpha} G),\mathbb{N},\beta)$. Therefore,
there exists a $*-$homomorphism
$$
\delta : (D[\varphi] \rtimes_{\alpha} G) \rtimes_\beta \mathbb{N} \longrightarrow  \mathbb{U}[\varphi]
$$
such that $\delta \circ i_D = \varpi$ and $\delta \circ i_N = \rho$.
\end{proof}

\section{OPEN QUESTIONS}

It is natural to ask if $\delta$ above is injective, which would imply that those two objects are isomorphic. We tried to find elements inside $\mathbb{U}[\varphi]$ which behave just like the images of $i_D$ and $i_N$, but it seems to be impossible since we would need that $su_{1_{\tau(h)}}s^*=u_{1_{\tau(\varphi(h))}}$. Perhaps imposing additional properties or requirements on the initial objects may be enough to guarantee this property.

Another point worth asking about is the $K$-theory of the structures presented here. Some developments have been made in the study of such structures constructed from simpler ones, as in \cite{Li4}.

%%%%%%%%%%%%%%%%%%%%%%%%%%%

\bibliography{ijmsample}
\bibliographystyle{ijmart}

\end{document}